\documentclass{amsart}
\usepackage[utf8]{inputenc}
\usepackage[utf8]{inputenc}
\usepackage[T1]{fontenc}
\usepackage[all]{xy}
\usepackage{lipsum}
\usepackage{url}
\usepackage{tikz}
\usepackage{stackrel}
\usepackage{color}
\def\red{\textcolor{black}}

\usepackage{amsmath,amsthm,amssymb}

\usepackage{xcolor}

\usepackage{graphicx}

\newtheorem{theorem}{Theorem}[section]

\newtheorem{example}[theorem]{Example}
\newtheorem{proposition}[theorem]{Proposition}
\theoremstyle{definition}

\newtheorem{remark}[theorem]{Remark}

\numberwithin{equation}{section}

\begin{document}

\vspace{0.5in}

\renewcommand{\bf}{\bfseries}
\renewcommand{\sc}{\scshape}
\vspace{0.5in}

\title[Motion planning algorithms]%
{Sequential collision-free optimal motion planning algorithms in punctured Euclidean spaces \\ }

\author{Cesar A. Ipanaque Zapata}
\address{Departamento de Matem\'{a}tica,UNIVERSIDADE DE S\~{A}O PAULO
INSTITUTO DE CI\^{E}NCIAS MATEM\'{A}TICAS E DE COMPUTA\c{C}\~{A}O -
USP , Avenida Trabalhador S\~{a}o-carlense, 400 - Centro CEP:
13566-590 - S\~{a}o Carlos - SP, Brasil}
\curraddr{Departamento de Matem\'{a}ticas, CENTRO DE INVESTIGACI\'{O}N Y DE ESTUDIOS AVANZADOS DEL I. P. N.
Av. Instituto Polit\'{e}cnico Nacional n\'{u}mero 2508,
San Pedro Zacatenco, Mexico City 07000, M\'{e}xico}
\email{cesarzapata@usp.br}
\thanks{The first author would like to thank grant\#2018/23678-6, S\~{a}o Paulo Research Foundation (FAPESP) for financial support.}

\author{Jes\'{u}s Gonz\'{a}lez}
\address{Departamento de Matem\'{a}ticas, CENTRO DE INVESTIGACI\'{O}N Y DE ESTUDIOS AVANZADOS DEL I. P. N.
Av. Instituto Polit\'{e}cnico Nacional n\'{u}mero 2508,
San Pedro Zacatenco, Mexico City 07000, M\'{e}xico}
\email{jesus@math.cinvestav.mx}



\begin{abstract} In robotics, a topological theory of motion planning was initiated by M. Farber. The multitasking motion planning problem is new and its theoretical part via topological complexity has hardly been developed, but the concrete implementations are still non-existent, and in fact this work takes the first step in this last direction (producing explicit algorithms.) We present optimal motion planning algorithms which can be used in designing practical systems controlling objects moving in Euclidean space without collisions between them and avoiding obstacles. Furthermore, we present the multitasking version of the algorithms.
\end{abstract}

\maketitle


\section{Introduction}

Robot motion planning problem usually ignores dynamics and other differential constraints and focuses primarily on the translations and rotations required to move the robot. Here, we will have in mind an infinitesimal mass particle as an object (e.g., infinitesimally small ball). 

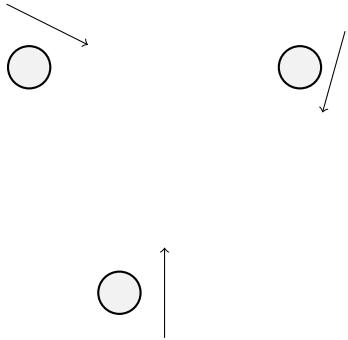
\begin{figure}[h]
 \centering
\begin{tikzpicture}[x=.6cm,y=.6cm]
\filldraw[color=black, fill=black!5, thick] (0.5,0) circle (8pt); 
\draw[->] (0,1.4) -- (1.8,0.5);
\filldraw[color=black, fill=black!5, thick] (6.5,0) circle (8pt); 
\draw[->] (7.5,0.8) -- (7,-1);
\filldraw[color=black, fill=black!5, thick] (2.5,-5) circle (8pt); 
\draw[->] (3.5,-6) -- (3.5,-4);
\end{tikzpicture}
\caption{Multi-robot system.}
 \label{multi-robot}
\end{figure}

Consider a multi-robot system consisting of $k$ distinguishable robots moving in Euclidean space $\mathbb{R}^d$ ($d\geq 2$) without collisions and avoiding $r$ stationary obstacles $(r\geq 0)$. In this work, we focus primarily on the translations required to move the robot. The associated state space or configuration space to this mechanical system is the classical ordered configuration space $F(\mathbb{R}^d-Q_r,k)$ of $k$ distinct points on punctured Euclidean space $\mathbb{R}^d-Q_r$ (see \cite{fadell1962configuration} for the notion of ordered configuration spaces). Here $Q_r=\{q_1,\ldots,q_r\}$ represents the set of the $r$ obstacles $q_j$. Explicitly, \[F(\mathbb{R}^d-Q_r,k)=\{(x_1,\ldots,x_k)\in (\mathbb{R}^d)^k\mid ~~x_i\neq x_j\text{ for } i\neq j \text{ and } x_i\neq q_j \text{ for any } i,j \},\] equipped with subspace topology of the Cartesian power $(\mathbb{R}^d)^k$. Note that the $i-th$ coordinate of a point $(x_1,\ldots,x_n)\in F(\mathbb{R}^d-Q_r,k)$ represents the state or position of the $i-th$ moving object, so that the condition $x_i\neq x_j$ reflects the collision-free requirement and the condition $x_i\neq q_j$ reflects the avoiding obstacle requirement.

\textit{The collision-free sequential robot motion planning problem} (a la Rudyak) consists in controlling simultaneously these $k$ robots without collisions between them and avoiding obstacles, where one is interested, in addition of initial-final states, in $n-2$ intermediate states of the robots. To solve this problem we need to find an \emph{$n$-th sequential collision-free optimal motion planning algorithm} on state space $F(\mathbb{R}^d-Q_r,k)$ (see Section 2). A central challenge of modern robotics (see, for example Latombe \cite{latombe2012robot} and LaValle \cite{lavalle2006planning}) consists of designing explicit and suitably optimal motion planners. This involves challenges in modeling planning problems, designing efficient algorithms, and developing robust implementations. These are exciting times to study planning algorithms and contribute to their development and use. 

Sequential collision-free optimal motion planning algorithms in Euclidean spaces without obstacles was given by the authors in \cite{zapata2019multitasking}. The purpose of the present work is to address the punctured case. 

In order to give sequential collision-free optimal motion planning algorithms, we need to know the smallest possible number of regions of continuity for any $n$-th sequential collision-free motion planning algorithm, that is, the value of $\text{TC}_n(F(\mathbb{R}^d-Q_r,k))$\footnote{In this work we use the non-reduced version of the $n-$sequential topological complexity, TC$_n$ (see Section \ref{preliminar}).}. This value was computed by Gonz{\'a}lez and Grant in \cite{gonzalez2015sequential}.

\begin{theorem}\emph{(\cite{gonzalez2015sequential})}\label{gongra}For $d,k,n\geq2$ and $r\geq 0$,
\[\text{TC}_n(F(\mathbb{R}^d-Q_r,k))=\left\{
  \begin{array}{ll}
     n(k-1)+1, & \hbox{if $r=0$ and $d$ is odd;}\\
     n(k-1), & \hbox{if $r=0$ and $d$ is even;}\\
     nk, & \hbox{if $r=1$ and $d$ is even;}\\
     nk+1, & \hbox{ otherwise.}
    \end{array}
\right.\]
\end{theorem}

For the experts, we can say that a higher optimal motion planning algorithm in $F(\mathbb{R}^d,k+1)$ induces a higher optimal motion planning algorithm in $F(\mathbb{R}^d-Q_1,k)$ (see Remark \ref{homotopy-invariance}) with $nk$ regions of continuity for $d\geq 2$. This is because $F(\mathbb{R}^d,k+1)$ and $F(\mathbb{R}^d-Q_1,k)$ are homotopy equivalent. Indeed we have $F(\mathbb{R}^d-Q_1,k)$ as a deformation retract of $F(\mathbb{R}^d,k+1)$. Thus, we focus in this work on the $r\geq 2$ case.

In this work we present a higher optimal motion planning algorithm in $F(\mathbb{R}^d-Q_r,k)$ with $nk+1$ regions of continuity. This algorithm works for any $d\geq 2$, $n\geq 2, r\geq 2$ and $k\geq 2$. Moreover, it gives an alternative proof in a constructive way to the inequality TC$_n(F(\mathbb{R}^d-Q_r,k))\leq nk+1$. This inequality was proved in \emph{\cite{gonzalez2015sequential}} using tools of homotopy theory.  

\section{Preliminary results}\label{preliminar}

The notion of $n$-th sequential or higher topological complexity was introduced by Rudyak in \cite{rudyak2010higher}, and further developed in~\cite{BGRT}. Here we follow \cite{zapata2019multitasking} to recall the basic definitions and properties.

\medskip
For a topological space $X$, let $PX$ denote the space of paths $\gamma:[0,1]\to X$, equipped with the compact-open topology. For $n\geq 2$, one has the evaluation fibration \begin{equation}\label{evaluation-fibration}
    e_n:PX\to X^n,~e_n(\gamma)=\left(\gamma(0),\red{\gamma\left(\dfrac{1}{n-1}\right)},\ldots,\red{\gamma\left(\dfrac{n-2}{n-1}\right)},\gamma(1)\right).
\end{equation} 
Recall, that an \textit{$n$-th sequential motion  planning  algorithm} is  a  section $s\colon X^n\to PX$ of  the  fibration  $e_n$, i.e.~ a (not necessarily continuous) map satisfying $e_n\circ s=id_{X^n}$. A continuous $n$-th sequential motion planning algorithm in $X$ exists if and only if the space $X$ is contractible. This fact gives, in a natural way, the definition of the following numerical invariant. The \textit{$n$-th sequential topological complexity} $\text{TC}_n(X)$ of a path-connected space $X$ is the Schwarz genus of the evaluation fibration~(\ref{evaluation-fibration}). In  other  words the $n$-th sequential topological complexity of $X$ is the smallest positive integer $\text{TC}_n(X)=k$ for which  the product $X^n$ is covered by $k$ open subsets $X^n=U_1\cup\cdots\cup U_k$ such that for any $i=1,2,\ldots,k$ there exists a continuous section $s_i:U_i\to PX$ of $e_n$ 
over $U_i$ (i.e., $e_n\circ s_i=i_{U_i}$), where $i_U:U\hookrightarrow X^n$ denotes the inclusion map. Any $n$-th sequential motion  planning  algorithm $s:=\{s_i:U_i\to PX\}_{i=1}^{k}$ is called \textit{optimal} if $k=\text{TC}_n(X)$.

\medskip One of the basic properties of TC$_n$ is its homotopy invariance, that is, if $X$ and $Y$ are homotopy equivalent then $\text{TC}_n(X)=\text{TC}_n(Y)$ for any $n\geq 2$. Furthermore, their $n$-th sequential motion  planning  algorithms are explicitly related.

\begin{remark}[Homotopy invariance]\label{homotopy-invariance}
Suppose $X$ \textit{dominates} $Y$, i.e., there are maps $f:X\to Y$ and $g:Y\to X$ such that $f\circ g\simeq id_Y$. Choose a homotopy $H:Y\times [0,1]\to Y$ with $H_0=id_Y$ and $H_1=f\circ g$. Let $\text{TC}_n(X)=k$ and let $s:=\{s_i:U_i\to PX\}_{i=1}^{k}$ be an $n$-th sequential motion  planning  algorithm to $X$ with $X^n=U_1\cup\cdots\cup U_k$ and $e_n\circ s_i=i_{U_i}$. For each $i=1,\ldots,k$ set $V_i:=(g\times\cdots\times g)^{-1}(U_i)\subseteq Y^n$ and define $\hat{s}_i:V_i\to PY$ by the formula
$$
  \hat{s}_i(y_1,\ldots,y_n)(t)= \begin{cases}
    H_{3(n-1)t}(y_1), & \hbox{$0\leq t\leq \frac{1}{3(n-1)}$;} \\
    f\left(s(g(y_1),\ldots,g(y_n))(3t-\frac{1}{n-1})\right), & \hbox{$\frac{1}{3(n-1)}\leq t\leq \frac{2}{3(n-1)}$;}\\
    H_{3-3(n-1)t}(y_2), & \hbox{$\frac{2}{3(n-1)}\leq t\leq \frac{1}{n-1}$;}\\
   H_{3(n-1)t-3}(y_2), & \hbox{$\frac{1}{n-1}\leq t\leq \frac{4}{3(n-1)}$;} \\
    f\left(s(g(y_1),\ldots,g(y_n))(3t-\frac{3}{n-1})\right), & \hbox{$\frac{4}{3(n-1)}\leq t\leq \frac{5}{3(n-1)}$;}\\
    H_{6-3(n-1)t}(y_3), & \hbox{$\frac{5}{3(n-1)}\leq t\leq \frac{2}{n-1}$;}\\  \quad\vdots\\
     H_{3(n-1)t-3(n-2)}(y_{n-1}), & \hbox{$\frac{n-2}{n-1}\leq t\leq \frac{3n-5}{3(n-1)}$;} \\
    f\left(s(g(y_{1}),\ldots,g(y_n))(3t-\frac{2n-3}{n-1})\right), & \hbox{$\frac{3n-5}{3(n-1)}\leq t\leq \frac{3n-4}{3(n-1)}$;}\\
   H_{3(n-1)-3(n-1)t}(y_{n}), & \hbox{$\frac{3n-4}{3(n-1)}\leq t\leq 1$.}
\end{cases}.
$$ One has $Y^n=V_1\cup\cdots\cup V_k$ and $e_n\circ \hat{s}_i=i_{V_i}$. Thus $\hat{s}:=\{\hat{s}_i:V_i\to PY\}_{i=1}^{k}$ is an $n$-th sequential motion planning algorithm to $Y$ and hence $\text{TC}_n(Y)\leq k=\text{TC}_n(X)$.

 In particular, if $X$ and $Y$ are homotopy equivalent we have $\text{TC}_n(X)=\text{TC}_n(Y)=k$. Furthermore, if $s:=\{s_i:U_i\to PX\}_{i=1}^{k}$ is  an optimal $n$-th sequential motion  planning  algorithm to $X$ then $\hat{s}:=\{\hat{s}_i:V_i\to PY\}_{i=1}^{k}$, as above, is an optimal $n$-th sequential motion  planning  algorithm to $Y$.
\end{remark}

\begin{remark}[Farber's TC and Rudyak's higher TC]
Note that $\text{TC}_2$ coincides with Farber`s topological complexity, which is defined in terms of motion planning algorithms for a robot moving between initial-final configurations~\cite{farber2003topological}. More general $\text{TC}_n$ is Rudyak`s higher topological complexity of motion planning problem, whose input requires, in addition of initial-final states, $n-2$ intermediate states of the robot. Similarly from \cite{zapata2019multitasking} we will use the expression ``motion planning algorithm'' as a substitute of ``$n$-th sequential motion planning algorithm for $n=2$''. 
\end{remark} 

Since~(\ref{evaluation-fibration}) is a fibration, the existence of a continuous motion planning algorithm on a subset $A$ of $X^n$ implies the existence of a corresponding continuous motion planning algorithm on any subset $B$ of $X^n$ deforming to $A$ within $X^n$. Such a fact is argued in \cite{zapata2019multitasking} in a constructive way, generalizing~\cite[Example 6.4]{farber2017configuration} (the latter given for $n=2$). This of course suits best our implementation-oriented objectives.

\begin{remark}[Constructing \red{motion planning algorithms} via deformations: higher case]\label{constructing-sections-via-deformations-higher-case}(\cite{zapata2019multitasking})
 Let \red{$s_A:A\to PX$ be a continuous motion planning algorithm defined on a subset $A$ of $X^n$.} Suppose a subset $B\subseteq X^n$ can be continuously deformed \red{within} $X^n$ into $A$. Choose a homotopy $H:B\times [0,1]\to X^n$ such that $H(b,0)=b$ and $H(b,1)\in A$ for any $b\in B$. Let $h_1,\ldots, h_n$ be the Cartesian components of $H$, $H=(h_1,\ldots,h_n)$. The formula
$$
 s_B(b)(\tau) = \begin{cases}
    h_1(b,3(n-1)\tau), & \hbox{$0\leq \tau\leq \frac{1}{3(n-1)}$;} \\
    s_A(H(b,1))(3\tau-\frac{1}{n-1}), & \hbox{$\frac{1}{3(n-1)}\leq \tau\leq \frac{2}{3(n-1)}$;}\\
    h_2(b,3-3(n-1)\tau), & \hbox{$\frac{2}{3(n-1)}\leq \tau\leq \frac{1}{n-1}$;}\\
    h_2(b,3(n-1)\tau-3), & \hbox{$\frac{1}{n-1}\leq \tau\leq \frac{4}{3(n-1)}$;} \\
    s_A(H(b,1))(3\tau-\frac{3}{n-1}), & \hbox{$\frac{4}{3(n-1)}\leq \tau\leq \frac{5}{3(n-1)}$;}\\
    h_3(b,6-3(n-1)\tau), & \hbox{$\frac{5}{3(n-1)}\leq \tau\leq \frac{2}{n-1}$;}\\  \quad\vdots\\
    h_{n-1}(b,3(n-1)\tau-3(n-2)), & \hbox{$\frac{n-2}{n-1}\leq \tau\leq \frac{3n-5}{3(n-1)}$;} \\
    s_A(H(b,1))(3\tau-\frac{2n-3}{n-1}), & \hbox{$\frac{3n-5}{3(n-1)}\leq \tau\leq \frac{3n-4}{3(n-1)}$;}\\
    h_n(b,3(n-1)-3(n-1)\tau), & \hbox{$\frac{3n-4}{3(n-1)}\leq \tau\leq 1$,}
\end{cases}.
$$
 defines a continuous section $s_B:B\to PX$ of~(\ref{evaluation-fibration}) over $B$. Hence, a deformation of $B$ into $A$ and a continuous motion planning algorithm defined on $A$ determine an explicit continuous motion planning algorithm defined on $B$.
\end{remark}

\subsection{Tame motion planning algorithms}
\medskip
Despite that the definition of $\text{TC}_n(X)$ deals with open subsets of $X^n$ admitting continuous sections of the evaluation fibration (\ref{evaluation-fibration}), for practical purposes, the construction of explicit $n$-th sequential motion planning algorithms is usually done by partitioning the whole space $X^n$ into pieces, over each of which~(\ref{evaluation-fibration}) a continuous section is set.  Since any such partition necessarily contains subsets which are not open (recall $X$ has been assumed to be path-connected), we need to be able to operate with subsets of $X^n$ of a more general nature.

Recall that a topological space $X$ is a \textit{Euclidean Neighbourhood Retract} (ENR) if it can be embedded into an Euclidean space $\mathbb{R}^d$ with an open neighbourhood $U$, $X\subset U\subset \mathbb{R}^d$, admitting a retraction $r:U\to X,$ $r\mid_U=id_X$.

\begin{example}
A subspace $X\subset \mathbb{R}^d$ is an ENR if and only if it is locally compact and locally contractible, see~\red{\cite[Chap.~4, Sect.~8]{dold2012lectures}}. This implies that all finite-dimensional polyhedra, smooth manifolds and semi-algebraic sets are ENRs.
\end{example}

Let $X$ be an ENR. We recall that a $n-th$ sequential motion planning algorithm $s:X^n\to PX$ is said to be \textit{tame} if $X^n$ splits as a pairwise disjoint union $X^n=F_1\cup\cdots\cup F_k$, where each $F_i$ is an ENR, and each restriction $s\mid_{F_i}:F_i\to PX$ is continuous. The subsets $F_i$ in such a decomposition are called \emph{domains of continuity} for $s$.

\begin{proposition}\emph{(\cite[Proposition 2.2]{rudyak2010higher})}\label{rudi}
For an ENR $X$, $\text{TC}_n(X)$ is the minimal number of domains of continuity $F_1,\ldots,F_k$ for tame $n$-th sequential motion planning algorithms $s:X^n\to PX$.
\end{proposition}

\begin{remark}
In the final paragraph of the introduction we noted that in this paper we construct optimal $n$-th sequential motion planners in $F(\mathbb{R}^d-Q_r,k)$. We can now be more precise: we actually construct $n$-th sequential tame motion planning algorithms with the advertized optimality property.
\end{remark}

\section{A tame motion planning algorithm in $F(\mathbb{R}^d-Q_r,k)$}

In this section we present a tame motion planning algorithm in $F(\mathbb{R}^d-Q_r,k)$ for $r\geq 2$. The algorithm with $2k+1$ regions of continuity works for any $d\geq 2,r\geq 2$ and $k\geq 2$;  this algorithm is optimal in the sense that it has the smallest possible number of regions of continuity. 

We think of $F(\mathbb{R},k+r)$ as a subspace of $F(\mathbb{R}^d,k+r)$ via the embedding $\mathbb{R}\hookrightarrow \mathbb{R}^d$, $x\mapsto (x,0,\ldots,0)$. By the $m-$homogeneous property of $\mathbb{R}^d$ we can suppose $Q_r=\{q_1,\ldots, q_r\}\subset \mathbb{R}$ with $q_1<q_2<\cdots <q_r$ and $\mid q_{i+1}-q_i\mid=1$. 

 Consider the first two standard basis elements $e_1=(1,0,\ldots,0)$ and $e_2=(0,1,0,\ldots,0)$ in $\mathbb{R}^d$ (we assume $d\geq2$). Denote by  $p:\mathbb{R}^d\to \mathbb{R},(x_1,\ldots,x_q)\mapsto x_1$ the projection onto the first coordinate. For  a  configuration $C\in F(\mathbb{R}^d,k+r)$,  where $C=(x_1,\ldots,x_{k+r})$ with $x_i\in\mathbb{R}^d, ~x_i\neq x_j$ for $i\neq j$, consider the set of projection points \[P(C) =\{p(x_1),\ldots,p(x_{k+r})\},\] $p(x_i)\in \mathbb{R}$, $i=1,\ldots,k+r$. The cardinality of this set will be denoted $\text{cp}(C)$. Note that $\text{cp}(C)$ can be any number $1, 2,\ldots,k+r$.

We recall that the configuration space $F(\mathbb{R}^d-Q_r,k)$ is the fiber of the Fadell-Neuwirth fibration $\pi_{k+r,r}:F(\mathbb{R}^d,k+r)\to F(\mathbb{R}^d,r),~(x_1,\ldots,x_{k+r})\mapsto (x_1,\ldots,x_{r})$. Indeed, the space $F(\mathbb{R}^d-Q_r,k)$ can be identify with the space \[\pi_{k+r,r}^{-1}(q_1,\ldots,q_r)=\{(q_1,\ldots,q_r,x_{r+1},\ldots,x_{r+k}):~~(x_{r+1},\ldots,x_{r+k})\in F(\mathbb{R}^d-Q_r,k)\}.\]

For our purposes, we recall the tame motion planning algorithm in $F(\mathbb{R}^d,k+r)$ given by the authors in \cite{zapata2019multitasking},  for any $d\geq 2$. This algorithm has domains of continuity $W_2,W_3,\ldots, W_{2k+2r}$, where \begin{equation*}
    W_l=\bigcup_{i+j=l}A_i\times A_j,
\end{equation*} and $A_i$ is the set of all configurations $C\in F(\mathbb{R}^d,k+r)$, with $\text{cp}(C)=i$.

We note that for each $i=1,\ldots,r-1$, the set $$A_i\cap F(\mathbb{R}^d-Q_r,k)=\varnothing,$$ because for any configuration $C\in F(\mathbb{R}^d-Q_r,k)=\pi_{k+r,r}^{-1}(q_1,\ldots,q_r)\subset F(\mathbb{R}^d,k+r)$, $\text{cp}(C)\geq r$. Then, for each $l=2,\ldots,2r-1$ we have $$W_l\cap \left(F(\mathbb{R}^d-Q_r,k)\times F(\mathbb{R}^d-Q_r,k) \right)=\varnothing.$$ Thus, in the rest of the paper we will consider $i\in\{r,r+1,\ldots,k+r\}$ and $l\in\{2r,2r+1,\ldots,2k+2r\}$. 

Set $A_i^\circ:=A_i\cap F(\mathbb{R}^d-Q_r,k)$ and $W_{l}^\circ=W_l\cap \left(F(\mathbb{R}^d-Q_r,k)\times F(\mathbb{R}^d-Q_r,k) \right)$. It is easy to see $W_{l}^\circ=\bigcup_{i+j=l}A_i^\circ\times A_j^\circ$. We recall from \cite{zapata2019multitasking}, the closure (relative to the topology on $F(\mathbb{R}^d,k+r)$) of each set $A_i$ is contained in the union of the sets $A_j$ with $j\leq i$:
\begin{equation*}
\overline{A_i}\subset \bigcup_{j\leq i}A_j.
\end{equation*} Hence, the closure (relative to the topology on $F(\mathbb{R}^d-Q_r,k)$) of each set $A_i^\circ$ is contained in the union of the sets $A_j^\circ$ with $j\leq i$:
\begin{equation}\label{cerraduras1}
\overline{A_i^\circ}\subset \bigcup_{j\leq i}A_j^\circ.
\end{equation}

The  sets $A_i^\circ=A_i\cap F(\mathbb{R}^d-Q_r,k)$ where $i=r, r+1,\ldots,k+r$,  are ENR, because they are semi-algebraic sets.

Next, we will construct a tame motion planning algorithm in $F(\mathbb{R}^d-Q_r,k)$ having $2k+1$ domains of continuity $W_{2r}^\circ,\ldots, W_{2k+2r}^\circ$. 

\subsection{Section over $F(\mathbb{R}-Q_r,k)\times F(\mathbb{R}-Q_r,k)$}

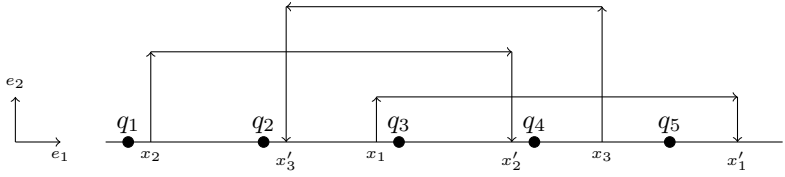
\begin{figure}[h]
 \centering
\begin{tikzpicture}[x=.6cm,y=.6cm]
\filldraw[black] (0.5,0) circle (2pt) node[above] {$q_1$};
\filldraw[black] (3.5,0) circle (2pt) node[above] {$q_2$};
\filldraw[black] (6.5,0) circle (2pt) node[above] {$q_3$};
\filldraw[black] (9.5,0) circle (2pt) node[above] {$q_4$};
\filldraw[black] (12.5,0) circle (2pt) node[above] {$q_5$};
\draw[->](-2,0)--(-1,0); \draw[->](6,0)--(6,1); \draw[->](1,0)--(1,2); \draw[->](11,0)--(11,3);
\draw[->](-2,0)--(-2,1); \draw[->](6,1)--(14,1); \draw[->](1,2)--(9,2); \draw[->](11,3)--(4,3);
\draw(0,0)--(15,0); \draw[->](14,1)--(14,0); \draw[->](9,2)--(9,0); \draw[->](4,3)--(4,0);
\node [below] at (-1,0) {\tiny$e_1$};
\node [above] at (-2,1) {\tiny$e_2$};
\node[below] at (1,0) {\tiny$x_2$};
\node[below] at (4,0) {\tiny$x_3^\prime$};
\node[below] at (6,0) {\tiny$x_1$};
\node[below] at (9,0) {\tiny$x_2^\prime$};
\node[below] at (11,0) {\tiny$x_3$};
\node[below] at (14,0) {\tiny$x_1^\prime$};
\end{tikzpicture}
\caption{Section over $F(\mathbb{R}-Q_r,k)\times F(\mathbb{R}-Q_r,k)$. Vertical arrows pointing upwards (downwards) describe
the first (last) third of the path $\Gamma^{C,C'},$ whereas horizontal arrows describe the middle third of $\Gamma^{C,C'}$.}
 \label{algorithm1}
\end{figure}

Given two configurations $C=(q_1,\ldots,q_r,x_{r+1},\ldots,x_{k+r})$ and $C^\prime=(q_1,\ldots,q_r,x^\prime_{r+1},\ldots,x^\prime_{k+r})$ in $F(\mathbb{R}-Q_r,k)$. Let $\Gamma^{C,C'}$ be the path in $F(\mathbb{R}^d-Q_r,k)$ from $C$ to $C^\prime$ depicted in Figure~\ref{algorithm1}.

Explicitly, $\Gamma^{C,C'}$ has components $(q_1,\ldots,q_r,\Gamma^{C,C'}_{r+1},\ldots,\Gamma^{C,C'}_{k+r})$ defined by 
\begin{equation}
    \label{Gamma1}\Gamma^{C,C'}_{i}(t)=\begin{cases}
    x_{i}+(3ti)e_2, & \hbox{for $0\leq t\leq \frac{1}{3}$;} \\
    x_{i}+ie_2+(3t-1)(x^\prime_{i}-x_{i}), & \hbox{for $\frac{1}{3}\leq t\leq \frac{2}{3}$;}\\
    x^\prime_{i}+i(3-3t)e_2, & \hbox{for $\frac{2}{3}\leq t\leq 1$}.
\end{cases}
\end{equation}
This yields a continuous motion planning algorithm $\Gamma:F(\mathbb{R}-Q_r,k)\times F(\mathbb{R}-Q_r,k)\to P F(\mathbb{R}^d-Q_r,k)$.

\subsection{The set $A_{k+r}^\circ$.}
If $C=(x_1,\ldots,x_r,x_{r+1},\ldots,x_{k+r})\in A_{k+r}$ then the map $\varphi:A_{k+r}\times [0,1]\to F(\mathbb{R}^d,k+r)$ given by the formula 
\begin{equation}
    \label{linear-transformation} \varphi_i(C,t)=x_i+t(p(x_i)-x_i),~i=1,\ldots,k+r
\end{equation}
defines a continuous deformation of $A_{k+r}$ onto $F(\mathbb{R},k)$ inside $F(\mathbb{R}^d,k)$ (see \cite{zapata2019multitasking}). We note that, if $C=(q_1,\ldots,q_r,x_{r+1},\ldots,x_{k+r})\in A_{k+r}^\circ= A_{k+r}\cap  F(\mathbb{R}^d-Q_r,k)$, $$\varphi_i(C,t)=q_i \text{ for } i=1,\ldots,r,$$ because $p(q_i)=q_i$. Thus, the restriction of $\varphi$ on $A_{k+r}\cap F(\mathbb{R}^d-Q_r,k)$ defines a continuous deformation of $A_{k+r}^\circ=A_{k+r}\cap F(\mathbb{R}^d-Q_r,k)$ onto $F(\mathbb{R}-Q_r,k)$ inside $F(\mathbb{R}^d-Q_r,k)$. As in Remark \ref{constructing-sections-via-deformations-higher-case}, this yields a continuous motion planning algorithm on $A_{k+r}^\circ\times A_{k+r}^\circ.$

\subsection{The sets $A_i^\circ$.} 
For a configuration $C\in A_i$, where $i\geq 2$, $C=(x_1,\ldots,x_{k+r})$ denote \[\epsilon (C):=\dfrac{1}{k+r}\min\{\mid p(x_r)-p(x_s)\mid:~ p(x_r)\neq p(x_s)\}. \] 

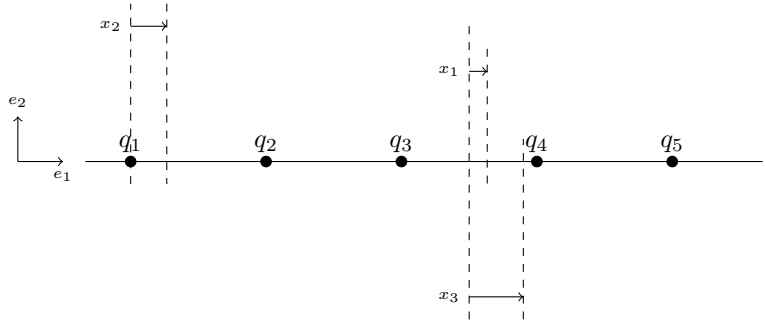
\begin{figure}[h]
 \centering
\begin{tikzpicture}[x=.6cm,y=.6cm]
\filldraw[black] (0.5,0) circle (2pt) node[above] {$q_1$};
\filldraw[black] (3.5,0) circle (2pt) node[above] {$q_2$};
\filldraw[black] (6.5,0) circle (2pt) node[above] {$q_3$};
\filldraw[black] (9.5,0) circle (2pt) node[above] {$q_4$};
\filldraw[black] (12.5,0) circle (2pt) node[above] {$q_5$};
\draw[dashed](0.5,-0.5)--(0.5,3.5); \draw[dashed](8,3)--(8,-3.5);
 \draw[->](8,2)--(8.4,2); \draw[->](0.5,3)--(1.3,3); \draw[->](8,-3)--(9.2,-3);
  \draw[dashed](8.4,2.5)--(8.4,-0.5); \draw[dashed](1.3,3.5)--(1.3,-0.5); \draw[dashed](9.2,-3.5)--(9.2,0.5);
 \draw[ ](-0.5,0)--(14.5,0);
\draw[->](-2,0)--(-2,1); \draw[->](-2,0)--(-1,0); 
\node [below] at (-1,0) {\tiny$e_1$};
\node [above] at (-2,1) {\tiny$e_2$};
\node[anchor=east] at (0.5,3) {\tiny$x_2$};
\node[anchor=east] at (8,2) {\tiny$x_1$};
\node[anchor=east] at (8,-3) {\tiny$x_3$};
\end{tikzpicture}
\caption{Desingularization deformation.}
 \label{algorithm2}
\end{figure}

For  $C\in A_i$ and $C\in F(\mathbb{R}^d-Q_r,k)$ (here we note $i\geq r\geq 2$), $C=(q_1,\ldots,q_r,x_{r+1},\ldots,x_{k+r})$, and $t\in[0,1]$, define 
\[D^i(C,t)=\begin{cases}(q_1,\ldots,q_r,z_{r+1}(C,t),\ldots,z_{k+r}(C,t)), & \mbox{if $r-1<i<k+r$;}\\
 C, & \mbox{if $i=k+r$,}
 \end{cases}\]
where $z_j(t)=x_j+t(j-1)\epsilon(C)e_1$ for $j=r+1,\ldots,k+r$. This  defines  a  continuous ``desingularization'' deformation $D^i:A_i^\circ\times [0,1]\to F(\mathbb{R}^d-Q_r,k)$ of $A_i^\circ$ into $A_{k+r}^\circ$ inside $F(\mathbb{R}^d-Q_r,k)$ despicted in Figure \ref{algorithm2}. 

Again, Remark \ref{constructing-sections-via-deformations-higher-case} yields a continuous motion planning algorithm on any subset $A_i^\circ\times A_j^\circ$ for $i,j\in\{r,r+1,\ldots,k+r\}$.

\subsection{Combining regions of continuity.} 
\begin{figure}[h]
 \centering
\begin{tikzpicture}[x=.6cm,y=.6cm]
\draw[->](-2,0)--(-1,0); 
\draw[->](-2,0)--(-2,1); 
\draw(0,0)--(17,0); 
\node [below] at (-1,0) {\tiny$e_1$};
\node [above] at (-2,1) {\tiny$e_2$};
\filldraw[black] (2,0) circle (2pt) node[above] {$q_1$};
\filldraw[black] (7,0) circle (2pt) node[above] {$q_2$};
 \draw[->](11,0)--(11.4,0);\draw[->](11.4,0)--(11.4,1); \draw[->](11.4,1)--(15,1); \draw[->](15,1)--(15,0); \draw[->](15,0)--(15,5); \draw[->](15,5)--(14,5); 
 \filldraw[black] (11,0) circle (2pt) node[above] {$x_1$};
 \filldraw[black] (14,5) circle (2pt) node[above] {$x_1^\prime$};
 \draw[->](4,-5)--(4.8,-5);\draw[->](4.8,-5)--(4.8,0); \draw[->](4.8,0)--(4.8,2); \draw[->](4.8,2)--(16,2); \draw[->](16,2)--(16,0); \draw[->](16,0)--(14,0); 
 \filldraw[black] (4,-5) circle (2pt) node[above] {$x_2$};
 \filldraw[black] (14,0) circle (2pt) node[above] {$x_2^\prime$};
 \draw[->](4,5)--(5.2,5);\draw[->](5.2,5)--(5.2,0); \draw[->](5.2,0)--(5.2,3); \draw[->](5.2,3)--(17,3); \draw[->](17,3)--(17,0); \draw[->](17,0)--(17,-5); \draw[->](17,-5)--(14,-5); 
 \filldraw[black] (4,5) circle (2pt) node[above] {$x_3$};
 \filldraw[black] (14,-5) circle (2pt) node[above] {$x_3$};
 \draw[dashed](4,-5)--(4,5);
 \draw[dashed](14,5)--(14,-5);
\end{tikzpicture}
\caption{The motion planning algorithm in $F(\mathbb{R}^d-Q_r,k)$.}
 \label{algorithm4}
\end{figure}
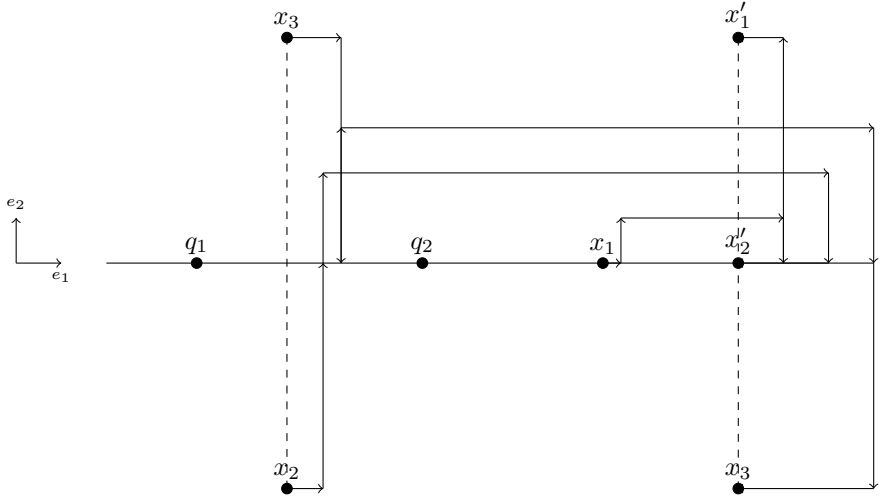

We have constructed continuous motion planning algorithms
\[\sigma_{i,j}\colon A_i^\circ\times A_j^\circ\to PF(\mathbb{R}^d-Q_r,k),\quad i,j=r, r+1,\ldots,k+r,\] by applying iteratively the construction in  Remark~\ref{constructing-sections-via-deformations-higher-case}. For $i,j\in\{r,r+1,\ldots,k+r\}$, the  sets $A_i^\circ\times A_j^\circ$ are pairwise disjoint ENR's covering $F(\mathbb{R}^d-Q_r,k)\times F(\mathbb{R}^d-Q_r,k)$. The resulting estimate $\text{TC}(F(\mathbb{R}^d-Q_r,k))\leq (k+1)^2$ is next improved by noticing that the sets $A_i^\circ\times A_j^\circ$  can be repacked into $2k+1$ pairwise disjoint ENR's each admitting its own continuous motion planning algorithm. Indeed,~(\ref{cerraduras1}) implies that $A_i^\circ\times A_j^\circ$ and $A_{i'}^\circ\times A_{j'}^\circ$ are ``topologically disjoint'' in the sense that $\overline{A_i^\circ\times A_j^\circ}\cap (A_{i^\prime}^\circ\times A_{j^\prime}^\circ)=\varnothing$, provided $i+j=i'+j'$ and $(i,j)\neq(i',j')$. Consequently, for $2r\leq\ell\leq 2k+2r,$ the motion planning algorithms $\sigma_{i,j}$ having $i+j=\ell$ determine a (well-defined) continuous motion planning algorithm on the ENR
\begin{equation*}
    W_{\ell}=\bigcup_{i+j=\ell}A_i^\circ\times A_j^\circ.
\end{equation*}
We have thus constructed a (global) tame motion planning algorithm in $F(\mathbb{R}^d-Q_r,k)$ having the $2k+1$ domains of continuity $W_{2r},W_{2r+1},\ldots, W_{2k+2r}$ (see Figure~\ref{algorithm4}).

\medskip
.

\begin{remark}
We note that the algorithm given here is not a restriction from the algorithm given in \cite{zapata2019multitasking}.
\end{remark}

\section{A higher tame motion planning algorithm in $F(\mathbb{R}^d-Q_r,k)$}\label{seccion4seqalgorithm}

In this section we present an optimal tame $n$-th sequential motion planning algorithm in $F(\mathbb{R}^d-Q_r,k)$, which generalizes in a natural way the algorithm presented in the previous section. As indicated in the introduction, the algorithm has $nk+1$ regions of continuity, works for any $d,k,n\geq 2$, and is optimal. The algorithm we present in this section can be used in designing practical systems controlling sequential motion of many objects moving in Euclidean space without collisions and avoiding obstacles.

\subsection{Section over $F(\mathbb{R}-Q_r,k)^n=F(\mathbb{R}-Q_r,k)\times\cdots\times F(\mathbb{R}-Q_r,k)$}
Recall we take the standard embedding $\mathbb{R}:=\{(x,0,\ldots,0)\in \mathbb{R}^d:~x\in \mathbb{R}\}$, so that $F(\mathbb{R}-Q_r,k)$ is naturally a subspace of $F(\mathbb{R}^d-Q_r,k)$. The motion planning algorithm $\Gamma:F(\mathbb{R}-Q_r,k)\times F(\mathbb{R}-Q_r,k)\to P F(\mathbb{R}^d-Q_r,k)$ given by~(\ref{Gamma1}) yields a continuous $n$-th motion planning algorithm \[\Gamma_n:F(\mathbb{R}-Q_r,k)\times\cdots\times F(\mathbb{R}-Q_r,k)\to P F(\mathbb{R}^d-Q_r,k)\] given by concatenation of paths \begin{equation}
    \Gamma_n(C_1,\ldots,C_n)=\Gamma(C_1,C_2)\ast\cdots\ast\Gamma(C_{n-1},C_{n}).
\end{equation}

\subsection{Motion planning algorithms $\sigma_{j_1,\ldots,j_n}$} We now go back to the notation introduced in the previous section where, for $r\leq i\leq k+r$, we constructed ENR's $A_i^\circ$ covering $F(\mathbb{R}^d-Q_r,k)$, as well as concatenated homotopies $A_i^\circ\times [0,1]\to F(\mathbb{R}^d-Q_r,k)$ deforming $A_i^\circ$ into $F(\mathbb{R}-Q_r,k)$. Together with the motion planning algorithm $\Gamma_n$, by Remark~\ref{constructing-sections-via-deformations-higher-case}, these deformations yield continuous $n$-th motion planning algorithms 
\[\sigma_{j_1,\ldots,j_n}:A_{j_1}^\circ\times\cdots\times A_{j_n}^\circ\to PF(\mathbb{R}^d-Q_r,k), \quad j_1,\ldots,j_n=r,r+1,\ldots,k+r.\] Indeed, the desingularization deformation $D^{j_1}\times\cdots\times D^{j_n}$ takes $A_{j_1}^\circ\times\cdots\times A_{j_n}^\circ$ into $(A_{k+r}^\circ)^n$; then we apply the deformation $\varphi\times\cdots\times\varphi \;\, (n-\text{times})$ which takes $(A_{k+r}^\circ)^n$ into $F(\mathbb{R}-Q_r,k)^n$; and finally we apply Remark \ref{constructing-sections-via-deformations-higher-case}. Let us emphasise that the above description of $\sigma_{j_1,\ldots,j_n}$ is fully implementable.

\subsection{Combining regions of continuity.} The  ENR's $A_{j_1}^\circ\times\cdots\times A_{j_n}^\circ$, $j_1,\ldots,j_n=r, r+1,\ldots,k+r$,  are mutually disjoint and cover the whole product $F(\mathbb{R}^d-Q_r,k)^n$. The resulting estimate $\text{TC}_n(F(\mathbb{R}^d-Q_r,k))\leq (k+1)^n$ coming from Proposition~\ref{rudi} and the motion planning algorithms $\sigma_{j_1,\ldots,j_n}$ is now improved by combining the domains of continuity to yield $nk+1$ covering ENR's $W_{\ell}$, $\ell=nr,nr+1,\ldots,nk+nr$, each admitting a continuous $n$-th motion planning algorithm. Explicitly, let \begin{equation}
    W_{\ell}=\bigcup_{j_1+\cdots+j_n=\ell}A_{j_1}^\circ\times\cdots\times A_{j_n}^\circ,
\end{equation}
where $\ell=nr,nr+1,\ldots,nk+nr$. By~(\ref{cerraduras1}), any two distinct $n-$tuples $(j_1,\ldots,j_n)$  and $(j^\prime_1,\ldots,j^\prime_n)$  with $j_1+\cdots+j_n=j^\prime_1+\cdots+j^\prime_n$ determine topologically disjoint sets $A_{j_1}^\circ\times\cdots\times A_{j_n}^\circ$ and $A_{j^\prime_1}^\circ\times\cdots\times A_{j^\prime_n}^\circ$ in $F(\mathbb{R}^d-Q_r,k)^n$, i.e. $\overline{A_{j_1}^\circ\times\cdots\times A_{j_n}^\circ}\cap (A_{j^\prime_1}^\circ\times\cdots\times A_{j^\prime_n}^\circ)=\varnothing$. Therefore the motion planning algorithms $\sigma_{j_1,\ldots,j_n}$ with $j_1+\cdots +j_n=\ell\hspace{.3mm}$ jointly  define  a  continuous motion planning algorithm on $W_{\ell}$. We have thus constructed a tame $n$-th sequential motion planning algorithm in $F(\mathbb{R}^d-Q_r,k)$ having $nk+1$ domains of continuity $W_{nr},W_{nr+1},\ldots, W_{nk+nr}$.

\bibliographystyle{plain}

\end{document}